\documentclass[reqno]{amsart}
\usepackage{bbold}
\usepackage{amsmath}
\usepackage{amssymb}
\usepackage{amsfonts}

\newcommand{\abs}[1]{\vert #1 \vert}

\newcommand{\norm}[1]{\left\Vert #1 \right\Vert}

\newcommand{\init}{\vert_{t = 0}}

\newcommand{\C}{\mathbb{C}}
\newcommand{\N}{\mathbb{N}}

\newcommand{\R}{\mathbb{R}}

\DeclareMathOperator{\im}{Im}
\DeclareMathOperator{\re}{Re}

\newtheorem{theorem}{Theorem}

\theoremstyle{definition}
\newtheorem{definition}{Definition}

\theoremstyle{remark}
\newtheorem{remark}{Remark}

\title[Ill-posedness for 1d Maxwell-Dirac]{Sharp ill-posedness for the Maxwell-Dirac equations in one space dimension}

\author[S.~Selberg]{Sigmund Selberg}
\author[A.~Tesfahun]{Achenef Tesfahun}

\address{Department of Mathematics, University of Bergen, PO Box 7803, 5020 Bergen, Norway}
\email{Sigmund.Selberg@uib.no}

\email{Achenef.Temesgen@uib.no}


\begin{document}

\begin{abstract}
The Maxwell-Dirac equations in one space dimension are proved to be well posed in the charge class, that is, with $L^2$ data for the spinor. We also prove that this result is sharp, in the sense that well-posedness fails for spinor data in $H^s$ with $s<0$, as well as in $L^p$ with $1 \le p < 2$. More precisely, we give an explicit example of such data for which no local solution can exist. Our proof of well-posedness applies to a class of systems which includes also the Dirac-Klein-Gordon system, but it does not require any null structure in the system.
\end{abstract}

\maketitle

\section{Introduction}

We consider the Maxwell-Dirac equations on the Minkowski space-time $\R^{1+1}$,
\begin{subequations}\label{MD}
\begin{align}
  \label{MDa}
  (-i\gamma^\mu \partial_\mu + M) \psi &= A_\mu \gamma^\mu \psi,
  \\
  \label{MDb}
  \square A_\mu &= - \overline{\psi} \gamma_\mu \psi,
\end{align}
\end{subequations}
with initial conditions at time $t=0$,
\begin{equation}\label{Data}
  \psi(0,x) = \psi_0(x), \qquad A_\mu(0,x) = a_\mu(x), \qquad \partial_t A_\mu(0,x) = b_\mu(x).
\end{equation}
The unknowns are the Dirac spinor field $\psi \colon \R^{1+1} \to \C^2$, regarded as a column vector, and the electromagnetic potential components $A_\mu \colon \R^{1+1} \to \R$, $\mu=0,1$. Here $\square = \partial^\mu \partial_\mu$ is the d'Alembertian, $M \in \R$ is a mass constant, and $\overline \psi = \psi^*\gamma^0$ with $\psi^*$ the complex conjugate transpose. The equations are written in covariant form on $\R^{1+1}$ with coordinates $x^\mu$ and metric $(g^{\mu\nu})=\mathrm{diag}(1,-1)$, where $x^0=t$ is time and $x^1=x$ is spatial position, and we write $\partial_\mu = \partial/\partial x^\mu$, so that $\partial_0=\partial_t$, $\partial_1=\partial_x$ and $\square = \partial_t^2-\partial_x^2$. The $2 \times 2$ Dirac matrices $\gamma^\mu$ should satisfy
\[
  \gamma^\mu \gamma^\nu + \gamma^\nu \gamma^\mu = 2g^{\mu\nu} I \qquad (g^{00}=1, g^{11}=-1, g^{01}=g^{10}=0)
\]
and
\[
  (\gamma^0)^* = \gamma^0, \qquad (\gamma^1)^* = -\gamma^1.
\]
We choose the representation
\[
  \gamma^0 =
  \begin{pmatrix}
    0 & 1  \\
    1 &0
  \end{pmatrix},
  \qquad
  \gamma^1=
  \begin{pmatrix}
    0 & -1  \\
    1 & 0
  \end{pmatrix}.
\]

The Maxwell-Dirac system describes the motion of an electron interacting with its self-induced electromagnetic field, and it is the fundamental PDE system in relativistic quantum electrodynamics.

A key fact about this system is that it enjoys a $\mathrm{U}(1)$ gauge freedom, and the particular form \eqref{MD} appears when the Lorenz gauge condition $\partial^\mu A_\mu = 0$ is chosen, that is,
\begin{equation}\label{LG}
  \partial_t A_0 = \partial_x A_1.
\end{equation}
Since the latter reduces to a constraint on the initial data, we do not include it in \eqref{MD}. A second, less obvious constraint on the data, arising from \eqref{MDb} and \eqref{LG} combined, is the Gauss law
\begin{equation}\label{Gauss}
  \partial_x E = \abs{\psi}^2,
\end{equation}
where
\[
  E := \partial_x A_0 -  \partial_t A_1
\]
is the electric field. If the constraints \eqref{LG} and \eqref{Gauss} are satisfied by the data at time $t=0$, then they will also be satisfied at all later times, for a sufficiently regular solution of \eqref{MD}. 

Another key feature of the Maxwell-Dirac system is the conservation of charge,
\begin{equation}\label{Charge}
  \int_{\R} \abs{\psi(x,t)}^2 \, dx = \int_{\R} \abs{\psi(x,0)}^2 \, dx,
\end{equation}
 for sufficiently regular solutions. For this reason, a solution for which the map $t \mapsto \psi(t,\cdot)$ is continuous into $L^2(\R)$ and satisfies \eqref{Charge}, will be referred to as a \emph{charge class solution}.
 
The final key property that we want to mention, is that in the massless case $M=0$, the system \eqref{MD} is invariant under the rescaling
\[
  \psi(t,x) \longrightarrow \lambda^{3/2} \psi(\lambda t, \lambda x), \qquad A_\mu(t,x) \longrightarrow \lambda A_\mu(\lambda t, \lambda x) \qquad (\lambda > 0).
\]
By the usual heuristics, this provides some information about possible obstructions to well-posedness in a given data space $X_0$. Specifically, if we send $\lambda$ to zero, then the existence time of the rescaled solution goes to infinity, and this is only reasonable if the $X_0$ norm of the rescaled data tends to zero, or at least stays bounded. A data space $X_0$ is called \emph{subcritical}, \emph{critical} or \emph{supercritical} according to whether the norm of the rescaled data tends to zero, remains constant or tends to infinity, respectively, as $\lambda$ tends to zero. In a supercritical data space $X_0$ one does not expect well-posedness to hold.

To see what this heuristic tells us in the case of the Maxwell-Dirac system, let us start with the $L^2$ based Sobolev spaces $H^s(\R)$. For data
\begin{equation}\label{Sobolev}
  (\psi_0,a_\mu,b_\mu) \in X_0 := H^s(\R) \times H^r(\R) \times H^{r-1}(\R),
 \end{equation}
 the critical regularity is seen to be $s = -1$ and $r=-1/2$ (for the homogeneous spaces), so based on scaling alone, one does not expect well-posedness if $s < -1$ or $r < -1/2$ (supercritical scaling). In fact, we shall see that there are far stronger restrictions on well-posedness than this, excluding the range $s < 0$. But before we get to this, let us mention some earlier results on well-posedness and ill-posedness of the Maxwell-Dirac system in one space dimension.
 
Chadam \cite{Chadam1973} proved local well-posedness of \eqref{MD} in the space \eqref{Sobolev} with $s=r=1$, and moreover using the conservation of charge he showed that the solution extends globally in time. Okamoto \cite{Okamoto2013} proved local well-posedness for $s > 0$, $r > 1/2$, $s \le r \le \min(s+1,2s+1/2)$  and $(s,r) \neq (1/2,3/2)$, thus barely failing to reach the point $(s,r) = (0,1/2)$. Moreover, he proved that for $s > 0$, the data-to-solution map fails to be $C^2$ if $r$ is outside the range specified above. In the massless case $M=0$, Okamoto also proved that the data-to-solution map fails to be continuous at the point $(s,r) = (0,1/2)$. This last result shows that, if one wants to prove well-posedness for $s=0$ (or below), the data for the electromagnetic potential $A_\mu$ cannot be taken in the Sobolev spaces. A result in this direction was obtained by Huh \cite{Huh2010} in the massless case $M=0$: Using the interesting fact that the system can then be explicitly integrated, he proved global existence of \eqref{MD} in the case $s=0$ with $a_\mu, b_\mu \in BC(\R)$, where $BC(\R)$ denotes the space of bounded and continuous functions. This is however not a well-posedness result, since $\partial_t A_\mu$ does not persist in the space $BC(\R)$. Global existence and uniqueness of weak solutions for $s=0$ with data $(a_\mu,b_\mu) \in L^\infty(\R) \times L^1(\R)$ was obtained by You and Zhang \cite{Zhang2014}, without the restriction to zero mass. But this is also not a well-posedness result, since continuity of the solution map is not proved, and it is also not proved that $\partial_t A_\mu$ persists in $L^1(\R)$.

Thus, no proper well-posedness result for the Cauchy problem \eqref{MD}, \eqref{Data} has been obtained previously in the charge class, that is for $\psi_0 \in L^2(\R)$ (see, however, Remark \ref{LWPremark} below). Here we prove such a result, with data for the potential $A_\mu$ taken in the following space.

\begin{definition}
Let $Y=Y(\R)$ be the space with norm $\norm{f}_{Y} = \norm{f}_{L^\infty(\R)} + \norm{f'}_{L^1(\R)}$.
\end{definition}

Thus, $Y$ is the space of absolutely continuous functions $f \colon \R \to \C$ with bounded variation (cf.\ Corollary 3.33 in \cite{Folland1999}), and $Y_{\mathrm{loc}}$ is the space of locally absolutely continuous functions.

Our first main result is then the following.

\begin{theorem}\label{Thm1}
The Cauchy problem \eqref{MD}, \eqref{Data} is globally well posed for initial data
\[
  (\psi_0,a_\mu,b_\mu) \in X_0 := L^2(\R) \times Y(\R) \times L^1(\R).
\]
That is, for any $T > 0$, the problem has a unique solution $(\psi,A_\mu)$ on $(-T,T) \times \R$, satisfying
\[
  (\psi,A_\mu,\partial_t A_\mu) \in C([-T,T];X_0).
\]
Moreover, the data-to-solution map is continuous from $X_0$ to $C([-T,T];X_0)$, and higher regularity persists. In particular, the solution is a limit in $C([-T,T];X_0)$ of smooth solutions.
\end{theorem}

\begin{remark}
The above data space has a subcritical scaling. In fact the scaling is the same as for the homogeneous version of \eqref{Sobolev} with $(s,r) = (0,1/2)$.
\end{remark}

 \begin{remark}
By persistence of higher regularity we mean that if, for some $N \in \mathbb N$, we have $\partial_x^j (\psi_0,a_\mu,b_\mu) \in X_0$ for $j \in \{0,\dots,N\}$, then it follows that $\partial_t^j \partial_x^k (\psi,A_\mu,\partial_t A_\mu) \in C([-T,T];X_0)$ for $j,k \in \{0,\dots,N\}$ with $j+k \le N$.
 \end{remark}
 
 \begin{remark}
So far, we did not take into account the data constraints \eqref{LG} and \eqref{Gauss}. Typically, these constraints are not compatible with the choice of data space for $\partial_t A_\mu$. Indeed, Okamoto \cite{Okamoto2013} observed that in Chadam's result \cite{Chadam1973}, the electric field $E = \partial_x A_0 - \partial_t A_1$ would initially belong to $L^2(\R)$, but this is not compatible with \eqref{Gauss}, which implies that $E(0,x) = c + \int_0^x \abs{\psi_0(x)}^2 \, dx$ is an increasing function in $x$. A similar incompatibility occurs in our Theorem \ref{Thm1}, since $E$ would belong to $L^1(\R)$ initially. However, these incompatibilities are easily resolved by using the finite speed of propagation and localising.
 \end{remark}
 
 \begin{remark}\label{LWPremark}
Another way of resolving the incompatibility issue discussed in the previous remark, is to use the constraints  \eqref{LG} and \eqref{Gauss} directly in the statement of the Cauchy problem. Then in \eqref{Data} one has the constraints
\[
  b_0 = \frac{d}{dx} a_1, \qquad b_1 = \frac{d}{dx} a_0 - E_0,
\]
where the initial value  $E_0$ of the electric field is required to satisfy the Gauss law \eqref{Gauss}. Then the initial data are $(\psi_0,a_0,a_1,E_0)$. Global well-posedness of \eqref{MD} with such data was proved by the first author in \cite{Selberg2018} with $\psi_0 \in L^2(\R)$ and $a_0,a_1,E_0 \in BC(\R)$.
\end{remark}

 Our next main result is that Theorem \ref{Thm1} is sharp. For this we take
 \begin{equation}\label{badspinor}
   \psi_0(x) = \chi_{[-1,1]}(x) \frac{1}{\abs{x}^{1/2}}  \begin{pmatrix} 1 \\ 1 \end{pmatrix} \qquad (x \in \R, \; x \neq 0)
 \end{equation}
 where $\chi_{[-1,1]}$ is the characteristic function of the interval $[-1,1]$. Then
 \[
   \psi_0 \in L^p(\R) \quad \text{for $1 \le p < 2$},
 \]
so by the dual of the Sobolev embedding $H^r(\R) \subset L^q(\R)$ for $2 \le q < \infty$ and $r= 1/2-1/q$, it follows that also
\[
  \psi_0 \in H^s(\R) \quad \text{for $s < 0$}.
\]
But clearly $\psi_0$ fails to belong to $L^2(\R)$.
 
 \begin{theorem}\label{Thm2}
The Cauchy problem \eqref{MD}, \eqref{Data} is ill posed in
\[
  (\psi_0,a_\mu,b_\mu) \in X_0 := H^s(\R) \times D_0 \times D_1 \quad \text{for $s < 0$},
\]
and in
\[
  (\psi_0,a_\mu,b_\mu) \in X_0 := L^p(\R) \times D_0 \times D_1 \quad \text{for $1 \le p < 2$},
\]
regardless of the choice of spaces $D_0, D_1 \subset \mathcal D'(\R)$. In fact, with $\psi_0$ as in \eqref{badspinor} and with $a_\mu=b_\mu=0$ for $\mu=0,1$, the problem has no local solution near the origin in $\R^{1+1}$ which is a distributional limit of charge solutions. \end{theorem}
 
In the next two sections we give the proofs of well-posedness and ill-posedness, respectively. In fact, our proof of well-posedness applies to a fairly general class of systems which includes not only the Maxwell-Dirac system (MD) but also the Dirac-Klein-Gordon system (DKG) as special cases.
 
 \section{Global well-posedness in the charge class of generic systems of MD/DKG type}
 
 Here we prove Theorem \ref{Thm1}. In fact, we prove it for a more general system of the form
 \begin{subequations}\label{Generic}
 \begin{align}
   \label{GenericA}
  (-i\gamma^\mu \partial_\mu + M)\psi &= \sum_{j=1}^N V_j \gamma^0 B_j \psi,
  \\
  \label{GenericB}
  (\square + m^2)V_j &= \psi^* C_j \psi,
\end{align}
\end{subequations}
with initial conditions
\begin{equation}\label{GenericData}
    \psi(0,x) = \psi_0(x), \qquad V_j(0,x) = v_j(x), \qquad \partial_t V_j(0,x) = w_j(x)
\end{equation}
and unknowns $\psi \colon \R^{1+1} \to \C^2$ and $V = (V_1,\dots,V_N) \colon \R^{1+1} \to \R^N$. Here $N \in \N$, $m, M \in \R$ are constants, and the $B_j$ and $C_j$ are constant $2 \times 2$ hermitian matrices. The assumption $C_j^* = C_j$ guarantees that $V_j$ stays real valued given that its data are real valued. From \eqref{GenericA} and $B_j^* = B_j$ it then follows that $j^\mu := \psi^* \gamma^0 \gamma^\mu \psi$ satisfies $\partial_\mu j^\mu = 0$, hence the conservation of charge \eqref{Charge} holds.

We will prove the following result, which contains Theorem \ref{Thm1} as a special case.

\begin{theorem}\label{Thm3}
If $m=0$, the Cauchy problem \eqref{Generic}, \eqref{GenericData} is globally well posed for initial data
\[
  (\psi_0,v,w) \in X_0 := L^2(\R;\C^2) \times Y(\R;\R^N) \times L^1(\R;\R^N).
\]
In general (that is, not assuming $m=0$), the same result holds for data
\[
  (\psi_0,v,w) \in X_{0,\mathrm{loc}} := L^2_{\mathrm{loc}}(\R;\C^2) \times Y_{\mathrm{loc}}(\R;\R^N) \times L^1_{\mathrm{loc}}(\R;\R^N).
\]
\end{theorem}

\begin{remark}
For $m=0$, the second statement in the theorem is a consequence of the first statement and finite speed of propagation.
\end{remark}

\begin{remark}
Since we apply a contraction argument, we get well-posedness in the strong sense, including existence, uniqueness, and smooth dependence on the data. Moreover, higher regularity persists, so smooth initial data give a smooth solution.
\end{remark}

\begin{remark}
By invariance of the system \eqref{Generic} under the reflection $(t,x,M,B_j) \to (-t,-x,-M,-B_j)$, it suffices to prove Theorem \ref{Thm3} for positive times.
\end{remark}

\begin{remark}
The system \eqref{Generic} includes as special cases not only the Maxwell-Dirac system \eqref{MD} but also the Dirac-Klein-Gordon system (DKG)
 \begin{align*}
  (-i\gamma^\mu \partial_\mu + M)\psi &= \phi \psi,
  \\
  (\square + m^2)\phi &= \overline \psi \psi,
\end{align*}
for which Bournaveas \cite{Bournaveas2000} proved global well-posedness in the charge class, improving the earlier $H^1$-result of Chadam \cite{Chadam1973}. The proof of Bournaveas relies crucially on a null structure in the DKG system, whereas our proof of Theorem \ref{Thm3} does not require any such structure (of course, the two results are not quite identical, since the choice of data spaces for $\phi$ and $\partial_t\phi$ differs). On the other hand, the null structure in DKG is certainly necessary if one wants to go below the charge, and in fact it is possible to go down to $\psi_0 \in H^s$ for $s > -1/2$, but not further; see \cite{Machihara2010, Machihara2016, Machihara2018}.
 \end{remark}
 
 The remainder of this section is devoted to the proof of Theorem \ref{Thm3}. For convenience we rewrite the system in terms of the Dirac matrices $\alpha=\gamma^0\gamma^1$ and $\beta = \gamma^0$:
  \begin{subequations}\label{Generic'}
 \begin{align}
   \label{Generic'A}
  (-i\partial_t - i \alpha \partial_x +  M\beta)\psi &= \sum V_j B_j \psi,
  \\
  \label{Generic'B}
  (\partial_t^2 - \partial_x^2 + m^2)V_j &= \psi^* C_j \psi.
\end{align}
\end{subequations}

\subsection{Preliminaries}

In preparation for the proof we recall some pertinent facts.

\subsubsection{Estimates for the Klein-Gordon and wave equations} For
\[
  \left( \partial_t^2 - \partial_x^2 + m^2 \right) u = F, \qquad (u,\partial_t u)\init = (f,g),
\]
we recall the solution formula (see \cite[Section 4.1.3]{PA02})
\begin{align*}
  u(t,x) &= \frac{f(x+t) + f(x-t)}{2}
  - \frac{m^2 t}{2} \int_{x-t}^{x+t} \frac{J_1\left(m\sqrt{t^2-(x-y)^2}\right)}{m\sqrt{t^2-(x-y)^2}} f(y) \, dy
  \\
  &\quad
  + \frac12 \int_{x-t}^{x+t} J_0\left(m\sqrt{t^2-(x-y)^2}\right)g(y) \, dy
  \\
  &\quad + \frac12 \int_0^t \int_{x-(t-s)}^{x+t-s} J_0\left(m\sqrt{(t-s)^2-(x-y)^2}\right)F(s,y) \, dy \, ds,
\end{align*}
where $J_0(x) = \sum_{n=0}^\infty \frac{(-1)^n}{n!n!} \left( \frac{x}{2} \right)^{2n}$ and $J_1(x) = \sum_{n=0}^\infty \frac{(-1)^n}{n!(n+1)!} \left( \frac{x}{2} \right)^{2n+1}$ are the Bessel functions of the first kind. It is well known that $J_0(x)$ and $J_1(x)$ are $O(1)$ (in fact they are $O(1/\sqrt{x})$) as $x \to \infty$, hence $J_0(x), x^{-1}J_1(x) \le C < \infty$ for all $x \ge 0$. Thus, for $t > 0$,
\begin{equation}\label{KleinGordonEstimate}
  \norm{u(t)}_{L^\infty} \le C (1+m^2t^2) \norm{f}_{L^\infty} + C\norm{g}_{L^1}
  + C \int_0^t \norm{F(s)}_{L^1} \, ds.
\end{equation}

For $m=0$ one recovers D'Alembert's formula for the wave equation,
\[
  u(t,x) = \frac{f(x+t) + f(x-t)}{2} + \frac12 \int_{x-t}^{x+t} g(y) \, dy
  + \frac12 \int_0^t \int_{x-(t-s)}^{x+t-s} F(s,y) \, dy \, ds,
\]
and \eqref{KleinGordonEstimate} holds with $C=1$. Moreover, differentiating one obtains
\[
  \norm{\partial_x u(t)}_{L^1}, \norm{\partial_t u(t)}_{L^1} \le \norm{f'}_{L^1} + \norm{g}_{L^1}
  + \int_0^t \norm{F(s)}_{L^1} \, ds
\]
for $t > 0$, hence also
\begin{equation}\label{WaveEstimate}
  \norm{u(t)}_{Y} + \norm{\partial_t u(t)}_{L^1}
  \le 3 \left( \norm{f}_{Y} + \norm{g}_{L^1} + \int_0^t \norm{F(s)}_{L^1} \, ds \right).
\end{equation}

\subsubsection{Energy inequality for the Dirac equation}

Consider the Dirac equation
\[
  (-i\partial_t - i\alpha\partial_x + M\beta)\psi = F, \qquad \psi\init = f.
\]
Applying $i\partial_t - i\alpha\partial_x + M\beta$ to both sides and using $\alpha\beta+\beta\alpha=0$, $\alpha^2 = I$ and $\beta^2 = I$, one obtains
\[
  (\partial_t^2-\partial_x^2+M^2)\psi = G,
  \qquad (\psi,\partial_t\psi)\init=(f,g),
\]
where $G = (i\partial_t - i\alpha\partial_x + M\beta)F$ and $g = iF(0)-\alpha f' - iM\beta f$. Thus, the Klein-Gordon solution formula from the previous subsection applies, so assuming for the moment that $f$ and $F$ are smooth and compactly supported, it follows that $\psi$ is smooth and that $\psi(t)$ is compactly supported for each $t$. Now premultiply the Dirac equation by $i\psi^*$, take real parts, and use $\alpha^*=\alpha$ and $\beta^*=\beta$, to get $\partial_t \rho + \partial_x j = 2\re (i\psi^* F)$, where $\rho = \psi^*\psi$ and $j = \psi^*\alpha\psi$. Integration in $x$ gives
\[
  \frac{d}{dt} \int \psi^*\psi \, dx
  =
  2\re\int i\psi^* F \, dx
  \le 2\norm{\psi(t)}_{L^2} \norm{F(t)}_{L^2},
\]
implying the \emph{energy inequality},
\begin{equation}\label{EnergyInequality}
  \norm{\psi(t)}_{L^2} \le \norm{f}_{L^2} + \int_0^t \norm{F(s)}_{L^2} \, ds.
\end{equation}
By a density argument, the smoothness and support assumptions on $f$ and $F$ can now be removed, so that the inequality is valid for any $f \in L^2$ and $F \in L_t^1L_x^2$, in which case $\psi \in C(\R;L^2)$.

\subsection{Proof of Theorem \ref{Thm3}}

Solving for the potentials, we first prove local well-posedness for the non-linear and non-local Dirac equation thus obtained, with a time of existence depending on the $X_0$ norm of the data $(\psi_0,v,w)$. To obtain local well-posedness of the full system \eqref{Generic'} we then show that $(v,\partial_t v)$ persists in $Y \times L^1$ (or its local version if $m \neq 0$). Moreover the $Y \times L^1$ norm is a priori bounded on any finite time interval, and this together with the conservation of charge implies that the local result extends globally.

\subsubsection{Step 1: Local well-posedness for a non-linear and non-local Dirac equation}

Fix the data $(\psi_0,v,w) \in X_0$. Solving for the $V_j$ in \eqref{Generic'}, we obtain
\begin{equation}\label{NonlinearDirac}
  (-i\partial_t - i\alpha\partial_x + M\beta)\psi = \sum \mathfrak V_j[\psi] B_j \psi,
  \qquad \psi(0,x)=\psi_0(x),
\end{equation}
where the operators $\psi \mapsto \mathfrak V_j[\psi]$ are given by
\begin{align*}
  \mathfrak V_j[\psi](t,x) &= \frac{v_j(x+t) + v_j(x-t)}{2}
  - \frac{m^2 t}{2} \int_{x-t}^{x+t} \frac{J_1\left(m\sqrt{t^2-(x-y)^2}\right)}{m\sqrt{t^2-(x-y)^2}} v_j(y) \, dy
  \\
  &\quad
  + \frac12 \int_{x-t}^{x+t} J_0\left(m\sqrt{t^2-(x-y)^2}\right)w_j(y) \, dy
  \\
  &\quad + \frac12 \int_0^t \int_{x-(t-s)}^{x+t-s} J_0\left(m\sqrt{(t-s)^2-(x-y)^2}\right)(\psi^* C_j \psi)(s,y) \, dy \, ds.
\end{align*}
From \eqref{KleinGordonEstimate} we see that for any $T > 0$,
\begin{align*}
  \norm{\mathfrak V[\psi]}_{C_T L^\infty}
  &\le C (1+m^2T^2) \norm{v}_{L^\infty} + C\norm{w}_{L^1}
  + C T \norm{\psi}_{C_T L^2}^2,
  \\
  \norm{\mathfrak V[\psi] - \mathfrak V[\psi']}_{C_T L^\infty}
  &\le
  C T \left( \norm{\psi}_{C_T L^2} + \norm{\psi'}_{C_T L^2} \right) \norm{\psi-\psi'}_{C_T L^2},
\end{align*}
where $C_T L^p = C([0,T];L^p(\R))$ with norm $\norm{u}_{C_T L^p} = \sup_{t \in [0,T]} \norm{u(t)}_{L^p}$. From these estimates and the energy inequality \eqref{EnergyInequality}, we now see that for a pair of equations in iterative form,
\begin{alignat*}{2}
  (-i\partial_t - i\alpha\partial_x + M\beta)\psi &= \sum \mathfrak V_j[\psi'] B_j \psi',&
  \qquad \psi(0,x) &= \psi_0(x),
  \\
  (-i\partial_t - i\alpha\partial_x + M\beta)\Psi &= \sum \mathfrak V_j[\Psi'] B_j \Psi',&
  \qquad \Psi(0,x) &= \psi_0(x),
\end{alignat*}
where $\psi',\Psi' \in C_TL^2$ (the previous iterates) are given, we get the estimates:
\begin{align*}
  \norm{\psi}_{C_T L^2}
  &\le
  \norm{\psi_0}_{L^2}
  +
  T \norm{\mathfrak V[\psi']}_{C_T L^\infty} \norm{\psi'}_{C_T L^2}
  \\
  &\le
  \norm{\psi_0}_{L^2}
  +
  C T (1+m^2T^2) \norm{(v,w)}_{L^\infty \times L^1} \norm{\psi'}_{C_T L^2}
  + CT^2 \norm{\psi'}_{C_T L^2}^3
\end{align*}
and
\begin{align*}
  \norm{\psi-\Psi}_{C_T L^2}
  &\le
  CT \norm{\mathfrak V[\psi'] - \mathfrak V[\Psi']}_{C_T L^\infty} \norm{\psi'}_{C_T L^2}
  \\
  &\quad
  +
  CT \norm{\mathfrak V[\Psi']}_{C_T L^\infty} \norm{\psi'-\Psi'}_{C_T L^2}
  \\
  &\le
  CT^2 \left( \norm{\psi'}_{C_T L^2} + \norm{\Psi'}_{C_T L^2} \right)^2 \norm{\psi'-\Psi'}_{C_T L^2}
  \\
  &\quad
  +
  CT (1+m^2T^2) \norm{(v,w)}_{L^\infty \times L^1} \norm{\psi'-\Psi'}_{C_T L^2},
\end{align*}
where $C$ changes from line to line and depends also on the matrices $B_j$. It now follows by a standard iteration argument that we have local well-posedness for \eqref{NonlinearDirac}, and for any $R > 0$ we have a time of existence $T = T(R) > 0$ for data with $\norm{\psi_0}_{L^2} + \norm{v}_{L^\infty} + \norm{w}_{L^1} \le R$. Moreover, conservation of charge holds, since this is true for smooth solutions with compactly supported data, and the solutions we obtain are limits in $C_TL^2$ of such solutions.

\subsubsection{Step 2: Persistence of $(V,\partial_t V)(t)$ in $Y \times L^1$ and global existence}

First take $m=0$. Then by \eqref{WaveEstimate} and conservation of charge,
\begin{align*}
  \norm{V(t)}_{Y} + \norm{\partial_t V(t)}_{L^1}
  &\le C\left( \norm{v}_{Y} + \norm{w}_{L^1} + \int_0^t \norm{\psi(s)}_{L^2}^2 \, ds \right)
  \\
  &\le C\left( \norm{v}_{Y} + \norm{w}_{L^1} + t \norm{\psi_0}_{L^2}^2 \right) = O(1+t),
\end{align*}
hence the local result extends globally.

Now consider $m \neq 0$. Then we are only claiming well-posedness in $X_{0,\mathrm{loc}}$, hence by finite speed of propagation we may assume that the data $(\psi_0,v,w) \in X_0$ are compactly supported, say in the interval $[-a,a]$. Then $(\psi,V,\partial_t V)(t)$ is supported in $[-a-t,a+t]$ for $t > 0$. Temporarily writing the equation \eqref{Generic'B} for $V_j$ as
\[
  (\partial_t^2 - \partial_x^2)V_j = -m^2 V_j + \psi^* C_j \psi,
\]
we apply \eqref{WaveEstimate} and obtain
\[
  \norm{V(t)}_{Y} + \norm{\partial_t V(t)}_{L^1}
  = O(1+t)
  + C m^2 t 2(a+t) \norm{V}_{L^\infty([0,t]\times\R)}.
\]
To control the last term we note that \eqref{KleinGordonEstimate} implies $\norm{V(t)}_{L^\infty} = O(1+t^2)$, hence $\norm{V(t)}_{Y} + \norm{\partial_t V(t)}_{L^1} = O(1 + t^4)$. This concludes the proof of Theorem \ref{Thm3}.
 
 \section{Ill-posedness of Maxwell-Dirac below charge}
 
In terms of the components of $\psi = (u,v)^\intercal$ and setting $A_+ = A_0+A_1$ and $A_-=A_0-A_1$, the system \eqref{MD} becomes
\begin{subequations}\label{MD'}
\begin{align}
  \label{MD'a}
(\partial_t + \partial_x)u &= -iMv + i A_+ u,
  \\
  \label{MD'b}
  (\partial_t - \partial_x)v &= -iMu + i A_- v,
  \\
  \label{MD'c}
  (\partial_t^2-\partial_x^2) A_+ &= - 2\abs{v}^2,
  \\
  \label{MD'd}
  (\partial_t^2-\partial_x^2) A_- &= - 2\abs{u}^2.
\end{align}
\end{subequations}
We take initial data
\begin{equation}\label{Data'}
  u(0,x) = f(x), \qquad v(0,x) = g(x), \qquad A_\pm(0,x) = \partial_t A_\pm(0,x) = 0.
\end{equation}
Then with
\begin{equation}\label{fg}
  f(x) = g(x) = \chi_{[-1,1]} (x)\frac{1}{\abs{x}^{1/2}},
\end{equation}
which belongs to $L^p(\R)$, $1 \le p < 2$, and to $H^s(\R)$, $s < 0$, we show ill-posedness by non-existence. More precisely, approximating with data
\begin{equation}\label{fgepsilon}
  f_\varepsilon(x) = g_\varepsilon(x) = \chi_{[-1,1]} (x)\frac{1}{(\varepsilon + \abs{x})^{1/2}} \quad \text{for $\varepsilon > 0$},
\end{equation}
and denoting by $(u_\varepsilon,v_\varepsilon,A_{+,\varepsilon},A_{-,\varepsilon})$ the corresponding charge solution (which exists globally by Theorem \ref{Thm1}), we show that $A_{+,\varepsilon}$ fails to have a limit in the sense of distributions as $\varepsilon \to 0$, in the region $t > \abs{x}$.

By the finite speed of propagation we may remove the characteristic function $\chi_{[-1,1]}(x)$ in the above data. Indeed, this does not affect the solution in the region $\abs{x} + \abs{t} \le 1$, and it suffices to prove the non-convergence in this region.

We first prove the massless case, $M = 0$. Then the system can be explicitly integrated, as observed in \cite{Huh2010}. The general case will then be handled by comparing the massive solution with the massless one.

\subsection{The massless case} Taking $M=0$, the system \eqref{MD'}, \eqref{Data'} is easily integrated.

First, integrating \eqref{MD'a} and \eqref{MD'b} along characteristics gives
\begin{align*}
  u(t,x) &= f(x-t) e^{i \phi_+(t,x)},
  \\
  v(t,x) &= g(x+t) e^{i \phi_-(t,x)},
\end{align*}
where
\begin{align*}
  \phi_+(t,x) &= \int_0^t A_+(\sigma,x-t+\sigma) \, d\sigma,
  \\
  \phi_-(t,x) &= \int_0^t A_-(\sigma,x+t-\sigma) \, d\sigma
\end{align*}
are real valued. Then, since $\abs{u(t,x)}^2 = \abs{f(x-t)}^2$ and $\abs{v(t,x)}^2 = \abs{g(x+t)}^2$, we can integrate \eqref{MD'c} and \eqref{MD'd} to get
\begin{align*}
  A_+(t,x) = - \int_0^t \int_{x-(t-s)}^{x+t-s} \abs{g(y+s)}^2 \, dy \, ds,
  \\
  A_-(t,x) = - \int_0^t \int_{x-(t-s)}^{x+t-s} \abs{f(y-s)}^2 \, dy \, ds.
\end{align*}

These formal computations are valid for well-posed solutions, in particular for the charge solutions $(u_\varepsilon,v_\varepsilon,A_{+,\varepsilon},A_{-,\varepsilon})$ with data as in \eqref{fgepsilon} (with the characteristic function removed, as remarked above), and one can now easily compute the complete solution. For our purposes, however, the following lower bound suffices: In the region $t > \abs{x}$,
\begin{align*}
  -A_{-,\varepsilon}(t,x) &=
  \int_0^t \int_{x-(t-s)}^{x+t-s} \frac{1}{\varepsilon+\abs{y-s}} \, dy \, ds
  \\
  &\ge
  \int_0^{\frac{x+t}{2}} \int_{s}^{x+t-s} \frac{1}{\varepsilon+y-s} \, dy \, ds
  \\
  &=
  \frac{x+t}{2} (-\log \varepsilon)
  + \frac12 (\varepsilon+x+t) \left[ \log  (\varepsilon+x+t) - 1 \right]
  - \frac12 \varepsilon ( \log\varepsilon - 1 ).
\end{align*}
Now fix a non-negative test function $\theta \in C_c^\infty(\R^2)$ supported in the region $t > \abs{x}$. Then it follows that
\[
  - \int A_{-,\varepsilon}(t,x) \theta(t,x) \, dt \, dx
  \ge (- \log \varepsilon) \int \frac{t+x}{2} \theta(t,x) \, dt \, dx + R_\varepsilon,
\]
where $R_\varepsilon$ converges, by the dominated convergence theorem, to
\[
  R = \int \left( \frac12 (x+t) \left[ \log  (x+t) - 1 \right] \right) \theta(t,x) \, dt \, dx,
\]
as $\varepsilon \to 0$. We conclude that $A_{-,\varepsilon}$ cannot converge in the sense of distributions on the region $t > \abs{x}$, and this proves Theorem \ref{Thm2} in the case $M=0$.

\subsection{The massive case} In the case $M \in \R$, $M \neq 0$, it suffices to show the lower bound, uniformly in $\varepsilon > 0$,
\begin{equation}
  \label{KeyBound}
  \abs{u_\varepsilon(t,x)}^2 \ge \frac12 \abs{f_\varepsilon(x-t)}^2 \quad \text{for $0 < t \ll 1$ and $t < x < 1-t$},
\end{equation}
since then for $\abs{x} < t \ll 1$ we obtain
\begin{align*}
  - A_{-,\varepsilon}(t,x) &= \int_0^t \int_{x-(t-s)}^{x+t-s} \abs{u_\varepsilon(s,y)}^2 \, dy \, ds
  \\
  &\ge \int_0^{\frac{x+t}{2}} \int_{s}^{x+t-s} \abs{u_\varepsilon(s,y)}^2 \, dy \, ds
  \\
  &\ge \frac12 \int_0^{\frac{x+t}{2}} \int_{s}^{x+t-s} \frac{1}{\varepsilon+y-s} \, dy \, ds,
\end{align*}
hence the argument from the case $M=0$ goes through and proves Theorem \ref{Thm2}.

So it only remains to prove \eqref{KeyBound}. To this end, observe that \eqref{MD'a} and \eqref{MD'b} imply
\begin{align*}
  (\partial_t + \partial_x) \abs{u}^2 &= - 2M \im (u\overline v),
  \\
  (\partial_t - \partial_x) \abs{v}^2 &= 2M \im (u\overline v),
\end{align*}
which integrates to
\begin{align}
  \label{u2}
  \abs{u(t,x)}^2 &= \abs{f(x-t)}^2 - 2M \int_0^t \im (u\overline v) (\sigma,x-t+\sigma) \, d\sigma,
  \\
  \label{v2}
  \abs{v(t,x)}^2 &= \abs{g(x+t)}^2 + 2M \int_0^t \im (u\overline v) (\sigma,x+t-\sigma) \, d\sigma.
\end{align}

Now fix $\varepsilon > 0$ and define, for $\rho > 0$,
\[
  B_\rho(t) = \sup_{t+\rho \le x \le 1-t} \left( \abs{u_\varepsilon(t,x)}^2 + \abs{v_\varepsilon(t,x)}^2 \right).
\]
Note that this quantity is finite, since the solution is smooth in the region $x > t > 0$.

Applying \eqref{u2} and \eqref{v2} we then find
\[
  B_\rho(t) \le \frac{2}{\varepsilon+\rho} + 2\abs{M} \int_0^t B_\rho(\sigma) \, d\sigma
\]
so by Gr\"onwall's inequality,
\[
  B_\rho(t) \le \frac{2}{\varepsilon+\rho} e^{2\abs{M}t} \le \frac{4}{\varepsilon+\rho}
\]
if $0 < t < (2\abs{M})^{-1}$, which we assume from now on.

Applying \eqref{u2} again we now conclude that, for $t+\rho \le x \le 1-t$,
\begin{align*}
  \abs{u_\varepsilon(t,x)}^2 &\ge \frac{1}{\varepsilon + x-t} - \abs{M} \int_0^t B_\rho(\sigma) \, d\sigma
  \\
  &\ge \frac{1}{\varepsilon+x-t} - \abs{M} t \frac{4}{\varepsilon+\rho}.
\end{align*}
Choosing $\rho = x-t$ we obtain
\[
  \abs{u_\varepsilon(t,x)}^2 \ge \frac{1/2}{\varepsilon+x-t} \quad \text{for $0 < t < \frac{1}{8\abs{M}}$ and $t < x < 1-t$},
\]
proving \eqref{KeyBound}.

\bibliographystyle{amsplain}
\bibliography{database}

\end{document}